\documentclass[12pt]{amsart}
\usepackage{amsfonts,amssymb,amscd,amsmath,amsthm,enumerate,verbatim,calc}

\usepackage{times} 
\usepackage{mathptm}

\theoremstyle{plain}
\newtheorem{Theorem}{Theorem}[section]

\newtheorem{Lemma}[Theorem]{Lemma}
\newtheorem{Corollary}[Theorem]{Corollary} 
\newtheorem{Proposition}[Theorem]{Proposition}
\newtheorem*{Question}{Question}
\newtheorem*{Conjecture}{Conjecture}
\theoremstyle{definition}
\newtheorem{Definition}[Theorem]{Definition}

\newtheorem{Remark}[Theorem]{Remark}

\theoremstyle{remark}

\newtheorem*{chunk*}{}

\numberwithin{equation}{Theorem}

\newcommand{\height}{\operatorname{height}}
\newcommand{\mlabel}[1]%
  {\mbox{}\marginpar{\raggedleft\hspace{0pt}{\rm\ttfamily#1}}\label{#1}}

\newcommand{\into}{\operatorname{\hookrightarrow}}

\newcommand{\e}{\operatorname{e}}

\newcommand{\length}{\operatorname{\lambda}}

\newcommand{\Spec}{\operatorname{Spec}}
\newcommand{\Reg}{\operatorname{Reg}}
\newcommand{\Sing}{\operatorname{Sing}}
\newcommand{\Max}{\operatorname{Max}}
\newcommand{\fm}{{\mathfrak m}}
\newcommand{\fn}{{\mathfrak n}}

\newcommand{\ringR}{\text{$(R,\fm,k)$ }}
\newcommand{\ringS}{\text{$(S,\fn,k)$}}
\newcommand{\ringA}{\text{$(A,\fm)$}}

\newcounter{hours}\newcounter{minutes}

\newcommand{\excise}[1]{}

\begin{document}
\title[On the upper semi-continuity of the Hilbert-Kunz multiplicity]
{On the upper semi-continuity of the Hilbert-Kunz multiplicity}

\author[F.~Enescu]{Florian Enescu}
\author[K.~Shimomoto]{Kazuma Shimomoto}
\address{Department of Mathematics and Statistics, Georgia State University, Atlanta,
GA 30030 USA and The Institute of Mathematics of the Romanian Academy, Bucharest, Romania}
\email{fenescu@mathstat.gsu.edu}
\address{Department of Mathematics, University of Utah, Salt Lake City,
UT  84112 USA}
\email{shimomot@math.utah.edu}

\maketitle

\begin{abstract}
We show that the Hilbert-Kunz multiplicity of a $d$-dimensional nonregular complete intersection over $\overline{\mathbf{F}_p}$, $p>2$ prime, is bounded by below by the Hilbert-Kunz multiplicity of $\sum _{i=0}^{d}  x_i^2=0$, answering positively a conjecture of Watanabe and Yoshida in the case of complete intersections.
      
\end{abstract}

\bigskip
\section{Introduction}

Let $(R,\fm)$ be a local ring containing a field of positive characteristic $p>0$. If $I$ is an ideal in $R$, then $I^{[q]}=(i^q: i \in I)$, where $q=p^e$ is a power of the characteristic. Let $R^{\circ} = R \setminus \cup P$, where $P$ runs over the set of all minimal primes of $R$. An element $x$ is said to belong to the {\it tight closure} of the ideal $I$ if there exists $c \in R^{\circ}$ such that $cx^q \in I^{[q]}$ for all sufficiently large $q=p^e$. The tight closure of $I$ is denoted by $I^\ast$. By a ${\it parameter \ ideal}$ we mean an ideal generated by a full system of parameters in $R$. For an $\fm$-primary ideal $I$, one can consider the Hilbert-Samuel multiplicity and the Hilbert-Kunz multiplicity. A ring $R$ is called unmixed if 
${\rm dim} (R/Q) = {\rm dim} (R)$, for all associated primes $Q$ of $R$.

\begin{Definition}
 Let $I$ be an $\fm$-primary ideal in a $d$-dimensional local ring $(R,\fm)$.
In what follows $\length(-)$ denotes the length function.

{\it The Hilbert-Kunz multiplicity of $R$ at $I$} is defined by $\e _{HK} (I)= \e _{HK}(I,R): = \displaystyle\lim_{q \to \infty}  \frac{\length(R/I^{[q]})}{q^d}$.  Monsky has shown that this limit exists and  is positive. If $I =\fm$, then we call $\e_{HK} (\fm, R)$ the Hilbert-Kunz multiplicity of $R$ and denote it by $\e_{HK}(R)$.

{\it The Hilbert-Samuel multiplicity of $R$ at $I$} is defined by $\e (I)= \e (I,R) := \displaystyle\lim_{n \to \infty} d! \frac{\length(R/I^n)}{n^d}$. The limit exists and it is positive and similarly $\e (\fm, R)$ is simply denoted $\e(R)$ and called the Hilbert-Samuel multiplicity of $R$.

\end{Definition}

It is known that for parameter ideals $I$, one has $\e(I) = \e_{HK}(I)$. The following sequence of inequalities is also known to hold:
$${\rm max} \{ 1, \frac{1}{d!} \e (I) \} \leq \e_{HK} (I) \leq \e(I)$$
for every $\fm$-primary ideal $I$.

By a result of Watanabe and Yoshida \cite{WY1},
an unmixed local ring $R$ of characteristic $p>0$ is
regular if and only if the Hilbert-Kunz multiplicity,
\[
    \e_{HK}(R)= 1.
\]
 
A short proof of this was given by Huneke and Yao in~\cite{HY}.

In~\cite{BE}, Blickle and Enescu have started a first investigation of the number 
\[
    \epsilon_{HK}(d,p) = \inf\{\e_{HK}(R)-1 : \text{$R$ non--regular,
    unmixed, $\dim R = d$, $char R = p$} \}.
\]
by showing that $\epsilon_{HK}(d,p)$ is always
\emph{strictly} positive, i.e\ the Hilbert-Kunz multiplicity of a
non-regular ring of fixed dimension and characteristic cannot be
arbitrarily close to one. They have raised the natural question whether
$\epsilon_{HK}(d,p)$ is attained. And if this is the case, what is the
significance of such rings with minimal Hilbert-Kunz multiplicity?

In~\cite{WY2}, Watanabe and Yoshida have formulated the following conjecture

\begin{Conjecture}[Watanabe-Yoshida]
\label{conjecture}
Let $d \geq 2$ and $p \neq 2$ prime. Put

$$R_{p,d}: = k[[X_0,...,X_d]]/(X_0^2+ \cdots + X_d^2).$$

Let $\ringR$ be a $d$-dimensional unmixed local ring and let $k = \overline {\mathbf{F_p}}$. Then the following statements hold:

\item
$(1)$ If $R$ is not regular, then $\e_{HK}(R) \geq \e_{HK}(R_{p,d})$.

\item
$(2)$ If  $\e_{HK}(R) = \e_{HK}(R_{p,d})$, then the $\fm$-adic completion of $R$ is isomorphic to $R_{p,d}$ as local rings.
\end{Conjecture}

The case $d=2$ has been solved affirmatively (see ~\cite{WY1, BE}). The cases $d=3,4$ are more difficult and have been answered affirmatively by Watanabe and Yoshida, ~\cite{WY2}. The case $d=1$ is easy to interpret since $\e_{HK} (A) = \e (A)$.

In this paper we would like to prove part (1) of the Conjecture for complete intersections.

We would like to finish the introduction by mentioning two results that will be needed later.

\begin{Proposition}[Kunz, 3.2 in~\cite{K1} and 3.9 in~\cite{K2}]
\label{kunz}

Let $\ringR \to \ringS$ be a flat local homomorphism of Noetherian rings of characteristic $p$ such that $S/\fm S$ is regular.

\item
$(1)$ If $x$ is part of a system of parameters on $R$ then $\e_{HK} (R) \leq \e_{HK}(R/xR)$.

\item
$(2)$ $\e_{HK}(R) = \e_{HK}(S)$.
\end{Proposition}

We should note that Watanabe and Yoshida (\cite{WY1}) gave an alternate proof of (1) under the assumption that $x$ is nonzerodivisor on $R$.

An element $f \in A[[t]]$ over a local ring $\ringA$ is called a ${\it distinguished \ polynomial}$ if $f = a_o + a_1 t + \cdots + a_{n-1} t+ t^n$, for some integer $n$ and $a_i \in \fm, i \geq 0$.

In what follows we will need the following classical result:

\begin{Theorem} [Weierstrass Preparation Theorem,~\cite{G}] 
Let $\ringA$ be a complete local ring and let $B=A[[t]]$. If $f= \sum_{i=0}^{\infty} a_i t^i \in B$ and if there exists $n \in \mathbf{N}$ such that $a_i \in \fm $ for all $i < n$ and $a_n \notin \fm$, then $f = u f_o$ where $u$ is a unit in $B$ and $f_o$ is a distinguished polynomial of degree $n$. Also, $u$ and $f_o$ are uniquely determined by $f$.
\end{Theorem}

We would like to thank Paul C.~Roberts for valuable advice with regard to this paper. We are grateful to the referee for helpful comments that enhanced our exposition. In particular, Lemma~\ref{claim} was suggested by the referee.
Also, Ian Aberbach and C\u at\u alin Ciuperc\u a have informed us that they have obtained Theorem~\ref{main} independently. While their methods do not use the dense upper-semicontinuity of the Hilbert-Kunz multiplicity, they resemble ours in spirit.

\section{Dense upper-semicontinuity of the Hilbert-Kunz Multiplicity}

Let $R$ be an equidimensional ring of characteristic $p >0$ such that $R$ is finite over $R^p$, i.e. $R$ is $F$-finite. Kunz has shown that if $R$ is $F$-finite, then $R$ is excellent.

We would like to discuss here several aspects of the Hilbert-Kunz multiplicity.
 E.~Kunz has shown that the function $f_e : \Spec(R) \to \mathbf{Q}$ where $$f_e(P) = \length (R_P/ P^{[p^e]}R_P) / p^{e \height(P)}$$ is upper-semi continuous on $\Spec(R)$ (Corollary 3.4 in~\cite{K2}).

\begin{Definition}
Let $\e_{HK} : \Spec(R) \to \mathbf{R}$, defined by $$\e_{HK}(P) : = \e_{HK}(PR_P, R_P).$$ We caution the reader that, although one can talk about the Hilbert-Kunz multiplicity of an ideal primary to the maximal ideal in a local ring, the notation just introduced will always refer to the Hilbert-Kunz multiplicity of a local ring, $R_P$, at its maximal ideal. 
Clearly, $\e_{HK}(P) = \lim_{e \to \infty} f_e(P)$. 
\end{Definition}

\begin{Question}
Is $\e_{HK}$ an upper-semi continuous function on $\Spec(R)$?
\end{Question}

It is known that $\e_{HK}(P) \leq \e_{HK}(Q)$ if $P \subset Q$ are prime ideals in $R$ (Proposition 3.3 in~\cite{K2}). However, this does not immediately imply that $\e_{HK}$ is upper-semi continuous.

\begin{Definition}

Let $T$ be a topological space. A function $f : T \to \mathbf{R}$ is called dense upper semi-continuous if for every $x$ in $T$ one can find a dense subset $U$ of $T$ containing $x$ such that $f(y) \leq f(x)$ for every $y \in U$.

\end{Definition}

We  would like to introduce some more definitions before stating our next result. In what follows, by a variety, we always mean an irreducible, reduced scheme defined over an algebraically closed field. For a linear system $\Gamma$ (complete or not) on a variety $X$ we can define a rational map $\phi_{\Gamma} : X \dasharrow \mathbf{P}^N$ by sending $x \in X$ to $[s_o (x): \cdots : s_N(x)]$, where $s_i$ form a $K$-basis of the system. $\Gamma$ is said to be composed of a pencil if the image of this map is one dimensional.  

\begin{Lemma}[First Theorem of Bertini,~\cite{FOV}, Theorem 3.4.10 ]
Let $X$ be a variety over $K$ and let $\Gamma$ be linear system which is not composed of a pencil such that its base locus has codimension at least $2$. Then the generic member of $\Gamma$ is irreducible.

\end{Lemma}

\begin{Corollary}
Let $X$ be a n-dimensional variety over $K$. Then for every $x, y$ in $X$ there is an irreducible curve $C$ that passes through $x$ and $y$.
\end{Corollary}

\begin{proof}
If $X$ is a curve then there is nothing to prove. Assume that $\dim X \geq 2$.

Consider the linear system $\Gamma$ consisting of all the hyperplane sections that pass through $x$ and $y$. Then by Bertini there is an irreducible member $X_1 \in \Gamma$ such that $x, y \in X_1$. Take the reduced structure of $X_1$ so that it is a variety, denoted by $(X_1)_{red}$. Again apply Bertini to $(X_1)_{red}$ to get irreducible $X_2$ chosen from the linear system consisting of all the hyperplanes passing through $x, y$ in $(X_1)_{red}$. Keeping this procedure, we obtain the chain of closed subvarieties, say
$$X \supseteq (X_1)_{red} \supseteq \cdots \supseteq (X_{n-1})_{red}$$
such that $(X_{n-1})_{red}$ is one-dimensional, irreducible, and contains $x, y$. 

Hence $(X_{n-1})_{red}$ is our desired curve.
\end{proof}

\begin{Theorem}
Let $K$ be an uncountable algebraically closed field and $R$ a finitely generated $K$-algebra which is equi-dimensional. Let $\Sing(R) \subset \Max(R)$ be the singular locus. Then $\e_{HK} : \Max(R) \to \mathbf{R}$ is dense upper semi-continuous on each component of $\Max(R)$. In particular, $\e_{HK} : \Max(R) \to \mathbf{R}$ is dense upper semi-continuous on each irreducible component of $\Sing(R)$. 
\end{Theorem}

\begin{proof}
$R$ is an excellent ring and hence the regular locus of $R$ is open. 

The case when $R$ is a domain goes as follows: the regular locus is non-empty (the zero ideal is in it) and, for each $Q$ as in the hypothesis, one can take $\Lambda = \Reg(R) \cup \{ Q \}$. This is a dense set and $\e_{HK} (P) = 1 \leq \e_{HK} (Q) $ for every $P \in \Lambda$.

Now if $R$ is not a domain (and in particular if the regular locus happens to be empty) we have to argue differently:

We know that for every $e$ there exists an open set $Q \in \Lambda _e$ such that $f_e (P) \leq f_e (Q)$ for every $P \in \Lambda_e$ (see Corollary 3.4 in \cite{K1}).

We will take $\Lambda: = \cap _e \Lambda_e$ and show that $\Lambda$ is dense.

In the following, since we work on one component of Max(R), we may assume that Max(R) is irreducible but may possibly be non-reduced.

We need to show that, for every $x \in \Max(R)$ and every open set $x \in U$, $ U \cap \Lambda \neq \emptyset$ holds.  In other words, $U \cap _e \Lambda_e \neq \emptyset$. Then by Corollary applied to $\Max(R)_{red}$ there is an irreducible curve $C$ that passes through $x$ and $Q$ and set $\lambda_e = C \cap \Lambda_e$. Each $\lambda_e$ is open in $C$ and hence it is the complement of a finite set.

We have that $ (U \cap C)$ is an open set in $C$ containing $x$ and so $ (U \cap C) \cap \lambda_e \neq \emptyset$. Otherwise, $U \cap C$ is contained in the union of the complements of $\lambda_e$ which is a countable set. But $U \cap C$ is open in $C$ and hence it is definitely uncountable and therefore dense.

We have shown that $ (U \cap C) \cap \lambda_e \neq \emptyset$ which shows that $U \cap _e \Lambda_e \neq \emptyset$ must also be true.
The second statement follows from the similar argument by applying Bertini to irreducible component of $\Sing(R)_{red}$.
\end{proof}

Let $R_o =k[[x_1,...,x_n]]/(f)$ be an $n-1$-dimensional hypersurface ring and define an $n$-dimensional hypersurface ring $R= k[[x_1,...,x_n]][t]/(f+tg)$, where $g$ is a formal power series with $g \neq 0, g(0) =0, g  \notin k \cdot f$. Obviously, $t$ is a nonzerodivisor on $R$.

In this section, we would like to study the behavior of the Hilbert-Kunz multiplicity of the fibers of the natural homomorphism $k[t] \to R= k[[x_1,...,x_n]][t]/(f+tg)$. We will assume that $k$ is an uncountable algebraically closed and so all the maximal ideals of $k[t]$ are of the form $(t-\alpha)$, with $\alpha \in k$.
Let $t_{\alpha}= t-\alpha$. One can note that $R/(t_{\alpha})$ is a local ring isomorphic to $R_{\alpha} =k[[x_1,...,x_n]]/(f + \alpha g)$ which is a $n-1$-dimensional hypersurface. This makes $t_\alpha$ a nonzerodivisor on $R$, for every $\alpha \in k$. We would also like to note that every maximal ideal of $R$ is of the form $\fm_{\alpha} = (x_1,...,x_n,t-\alpha)$ with $\alpha \in k$.

\begin{Theorem}
\label{sc-hyp}
Assume that we are in the situation described above. 

One can find a dense subset $\Lambda \subset k$ such that, for every $\alpha \in \Lambda$, $$\e_{HK}( (R/t_{\alpha})_{\fm _{\alpha}})) = \e_{HK}(\frac{k[[x_1,...,x_d]]}{(f + \alpha g)}) \leq \e_{HK} ((R/tR)_{\fm _0})) = \e_{HK}( \frac{k[[x_1,..,x_d]]}{(f)}),$$

where $\fm_0 = (x_1,...,x_n,t)$.

\end{Theorem}

\begin{proof}
As remarked earlier, $R/t_{\alpha}R$ is already local with maximal ideal $\fm _{\alpha}$.

If $\ringA$ is a local ring of dimension $d$, the $\e_{HK} (A) = \lim _{q \to \infty} \length (A/ \fm ^{[q]}) / q^d$. Since $R/t_{\alpha}R$ and $R/tR$ have the same dimension, to prove the inequality in the statement we need to prove the inequality between the corresponding lengths.

Let us observe that, for every $\alpha$, $R / (\fm _{\alpha} ^{[q]} + t_{\alpha})R = R/ (x_1,...,x_n)^{[q]} \otimes _{k[t]} k[t]/(t_{\alpha})$. 

Moreover, let $A= R/ (x_1,...,x_n)^{[q]}$ and note that this is a finitely generated module over $k[t]$. So, if we localize at the multiplicative set $k[t] \setminus (t_\alpha)$ we get that $A_{(t_\alpha)}$ is a finitely generated module over $k[t]_{(t_\alpha)}$. Moreover, $A/(t_\alpha)$ is already local and we have that
$A/(t_\alpha) \simeq (A/(t_\alpha)) _{(t_\alpha)}$.

Since $k$ is algebraically closed, $\length (R / (\fm _{\alpha} ^{[q]} + t_{\alpha})R)$ equals the dimension of the $k$-vector space $R / (\fm _{\alpha} ^{[q]} + t_{\alpha})R = A/(t_\alpha)$. This, by NAK lemma, equals the minimal number of generators of $(R/ (x_1,...,x_n)^{[q]}) _{(t_{\alpha})} = A_{(t_\alpha)}$ over $k[t] _{(t_{\alpha})}$.

So, if we start with a set of minimal generators of $A _{(t)}$ over $k[t]_{(t)}$ we can find an open set $\Lambda_q$ in $k$, containing $0$, where we can extend these generators. 

Let $\Lambda = \cap _{q} \Lambda_q$. Since $k$ is uncountable and the complements of $\Lambda_q$ are all finite we see that $\Lambda$ must be an uncountable set and hence dense in $k$ in the Zariski topology.

For all $\alpha \in \Lambda$ we have that, for all $q$,

$$ \length (R / (\fm _{\alpha} ^{[q]} + t_{\alpha})R) \leq \length (R / (\fm _{0} ^{[q]} + t_{0})R), $$

and this gives the inequality that we want.
\end{proof}

We would like to close this section by discussing an example by Monsky that shows that one cannot hope to replace dense upper semi-continuity by upper semi-continuity in Theorem~\ref{sc-hyp}. 

First we would like to recall Monsky's example (\cite{M}):

\begin{Theorem}[Monsky]

Let $k$ be an algebraically closed field of characteristic $2$ and $R_\alpha=k[[x,y,z]]/(f+\alpha g)$, where $f = z^4+xyz^2+(x^3+y^3)z$, $g=x^2y^2$ and $0 \neq \alpha \in k$.

Then $\e_{HK} (R_\alpha) = 3+ 4^{-m_\alpha}$, where $m_\alpha$ is computed as follows. Write $\alpha = \beta ^2 + \beta$ with $\beta \in k$. 

\item 
$(1)$ If $\alpha$ is algebraic over $\mathbf{Z}/2\mathbf{Z}$, then $m_\alpha$ is the degree of $\beta$ over $\mathbf{Z}/2\mathbf{Z}$.

\item
$(2)$ If $\alpha$ is not algebraic over $\mathbf{Z}/2\mathbf{Z}$, then let $m_\alpha = \infty$.

\end{Theorem}

We would like to consider the case when $k$ is the algebraic closure of $(\mathbf{Z}/2\mathbf{Z}) (w)$, where $w$ is an indeterminate. Let $R = k[[x,y,z,t]]/(f + tg)$. We see that $R_\alpha = R/ (t-\alpha)$, where $\alpha \in k$.

We would like to show that $\e_{HK}$ is not necessarily upper semi-continuous in fibers over $k[t]$. More precisely, we will find $\alpha_0 \in k$ such that there exist no open subset $U$ in $k$ containing $\alpha_0$ such that  $\e_{HK} (R_\alpha) \leq \e_{HK}(R_{\alpha_0})$ for every $\alpha \in U$. If such $U$ exists, it would imply that $\e_{HK} (R_\alpha) > \e_{HK}(R_{\alpha_0})$ only for finitely many $\alpha$. However, if one takes $\alpha_0 = w$, we see that $\e_{HK}(R_{\alpha_0}) = 3$, because $w$ is not algebraic over $\mathbf{Z}/2\mathbf{Z}$. However, there are infinitely many elements $\alpha$ in $k$ that are algebraic over $\mathbf{Z}/2\mathbf{Z}$ and hence $\e_{HK}(R_\alpha) > 3$ for all these $\alpha$.

In conclusion, this example shows that if one wants to study the upper semi-continuity of the Hilbert-Kunz multiplicity of the fibers of $k[t] \to R$, a weaker notion of upper-semicontinuity must be considered. One example is our notion that replaces open sets by dense sets.

In what follows we will show how this notion can be exploited to prove a conjecture of Watanabe and Yoshida on the minimal Hilbert-Kunz multiplicity of non-regular rings.

\section{Minimal Hilbert-Kunz multiplicity: the hypersurface case}

\begin{Lemma}
\label{claim}
Let $k$ be a field such that $1/2 \in k$ and put $A=k[[x_1, ....,x_d]]$. Consider $B = A[[x_0]]$
and $F = x_0^2 + \cdots + x_d^2 +G$ with $G \in m_B^3$, where $m_B$ is the maximal ideal of $B$. 
Then there exist a unit $v_0$ in $B$, $a_0 \in (x_1,...,x_d)B$ and $G_1 \in (x_1,...,x_d)^3B$ such that 

$$F = v_0(x_0+a_0)^2+x_1^2+ \cdots + x_d^2 +G_1$$
\end{Lemma}

\begin{proof}
Write 
$$ G = \sum_{i=0}^{\infty} c_i x_0^i,$$ such that $c_i \in A$ and $c_0 \in m_A^3$, $c_1 \in m_A^2$ and $c_2 \in m_A$.

Let $v_0 = (1+c_2) + \sum_{i=1}^{\infty} c_{i+2}x_0^i$ and note that this is a unit in $B$.
Moreover,
$$F = v_0x_0^2+c_1x_0+c_0+x_1^2+ \cdots+x_d^2.$$

Now, let $a_0 = 2^{-1}v_0^{-1}c_1$ and $G_1=c_0-v_0a_0^2$ and note that the conclusion of the Lemma follows.

\end{proof}

\begin{Theorem}
\label{hypersurface}
For any $d$-dimensional singular hypersurface $k[[x_0,...,x_d]]/(f)$ over an uncountable algebraically closed field $k$ of characteristic different than $2$, we have that 
$$\e_{HK}(k[[x_0,...,x_d]]/(\sum_{i=0}^{d} x_i^2)) \leq \e_{HK}(R).$$
\end{Theorem}

\begin{proof}

We can assume that $f = \sum_{i=0}^{\infty} f_i$ where each $f_i$ is a homogeneous polynomial of degree $i$ and $f_0=f_1 =0$.

Since the characteristic of $k$ is different than $2$, we can make a change of variables to have that $f_2 = \sum_{i=0}^{l} x_i^2$ for some $-1 \leq l \leq d$ where $l =-1$ means that $f_2 =0$. 

Let us take $g_\alpha : = \alpha (x_{l+1}^2 + \cdots x_d ^2)$ with $\alpha \in k$. By Theorem~\ref{sc-hyp}, the Hilbert-Kunz multiplicity of $f$ is greater or equal than that of $F_\alpha = f+ g_\alpha$ for a dense set of $\alpha$'s in $k$. We can rescale our indeterminates and assume that $F_\alpha = x_o^2+ \cdots + x_d^2 + G$, where the $G$ contains only terms of degree greater or equal to $3$.

Apply Lemma~\ref{claim} to $F_\alpha$ and write $F_{\alpha} = v_0(x_0+a_0)^2+x_1^2+\cdots+x_d^2+G_1$, with $G_1$ an element of $(x_1,...,x_d)^3.$ We can continue now with $x_1^2+\cdots+x_d^2+G_1$ and by applying Lemma~\ref{claim} recursively we see that eventually we can write $F_\alpha = \sum_{i=0}^{d} v_i x_i^2$, where $v_i$ are all units, after a suitable change of variables.

Since we are working over an algebraically closed field of characteristic different than $2$, we can find $w_i$ units in $k[[x_0,...,x_d]]$ such that $w_i^2 = v_i$ (see Lemma~\ref{powers}). This allows us to transform $F_\alpha$ isomorphically into $\sum_{i=0}^d x_i^2$.

In conclusion, we get that
$$\e_{HK}(k[[x_0,...,x_d]]/(\sum_{i=0}^{d} x_i^2)) \leq \e_{HK}(R).$$
\end{proof}

\begin{Lemma}
\label{powers}
If $A$ is a ring such that $f = \sum u_i x^i$ is a formal power series in $A[[x]]$ and $u_o$ is a unit in $A$ that admits a square root in $A$ and $1/2 \in A$, we can find $g \in A[[x]]$ such that $g^2 =f$. In particular, if $f \in k[[x_0,...,x_d]]$ is a unit and $k$ is algebraically closed of characteristic different than $2$, then there exists $g \in k[[x_0,...,x_d]]$ such that $g^2 =f$.
\end{Lemma}

\begin{proof}
The first statement amounts to solving a system of equations where the unknowns are the coefficients of $g$. 

The second statement reduces to the first, by thinking of $f \in A[[x_d]]$ where $A=k[[x_0,...,x_{d-1}]]$. First, we apply induction on $d$: since $f$ is a unit, by induction we see that its constant term (when thinking of it as a power series in $x_d$ only) has a square root in $A=k[[x_0,...,x_{d-1}]]$. Applying the first statement now, we can find a power series $g \in A[[x_d]]=k[[x_0,...,x_d]]$ such that $g^2 =f$.

\end{proof}

Using an argument similar to the one in the proof of Theorem~\ref{hypersurface}, one can show the following:

\begin{Theorem} Let $(R,\fm,k)$ be a d-dimensional singular hypersurface complete local ring of characteristic $p>0$ and $p \ne 2,3$. Then one of the following is true.

\item
$(1)$ $R \cong k[[x_{0},...,x_{d}]]/(\sum_{i=0}^{d} x^{2}_{i})$, or   
\item
$(2)$ $\e_{HK}(R) \ge \e_{HK}(k[[x_{0},...,x_{d}]]/(x^{2}_{0}+\cdots+x^{2}_{d-1}+x^{3}_{d})).$
\end{Theorem}

\begin{proof}
Suppose that $R$ is defined by some $f \in k[[x_{0},...,x_{d}]]$. 

Assume $(1)$ is not the case. Then as in the proof of Theorem3.1, we can make change of variables to have that $f_{2}=\sum_{i=0}^{l}x^{2}_{i}$ for the homogeneous decomposition $f=\sum_{i=0}^{\infty}f_{i}$ of $f$. Since $(i)$ is not the case, we have that $l<d$. 

Let us take $g_{\alpha}:=\alpha (x^{2}_{l+1}+\cdots+x^{2}_{d-1}+x^{3}_{d})$ with $\alpha\in k$. Then $F_{\alpha}:=f+g_{\alpha}$ is of the form $x^{2}_{0}+\cdots+x^2_l+ \alpha x^{2}_{l+1}+\cdots+\alpha x^{3}_{d}+G$ for $\alpha \ne 0$, where $G$ contains only terms of degree greater than 2. 

Now we can keep track of the proof in Theorem~3.1 without any change to have that $F_{\alpha}=v_{o}x^{2}_{0}+\cdots+v_{d-1}x^{2}_{d-1}+v_{d}x^{3}_{d}$, where $v_{i}$ are all units. Since we can assume that $k$ is an algebraically closed field, and the characteristic of $k$ is different than 2 and 3, we can apply Lemma~3.2 to solve the system of equations in $w_{i}$; $w^{2}_{0}=v_{0}$,...,$w^{2}_{d-1}=v_{d-1}$,and $w^{3}_{d}=v_{d}$ (This is where $p \ne 3$ is used.). Therefore $F_{\alpha}$ can be transformed isomorphically into $x^{2}_{0}+\cdots+x^{2}_{d-1}+x^{3}_{d}$.

By dense upper semi-continuity, we get that
$$\e_{HK}(R) \ge \e_{HK}(k[[x_{0},...,x_{d}]]/(x^{2}_{0}+\cdots+x^{2}_{d-1}+x^{3}_{d})).$$
\end{proof}

Much has been learned about the Hilbert-Kunz multiplicity in Noetherian rings by comparing it to the more classical notion of Hilbert-Samuel multiplicity. It is true that in many instances the behavior of these two multiplicities is similar to each other. 

A natural way of approaching the conjecture of Watanabe and Yoshida is to show that for any equidimensional local ring $R$ there is a hypersurface $S$ of same dimension such that $\e_{HK} (S) \leq \e_{HK}(R)$. A well-known result on the Hilbert-Samuel multiplicity says that for every ring $R$ of dimension $d$ one can naturally construct, through Noether normalization, a $d$-dimensional hypersurface $S$ such that $\e(R) = \e(S)$. In this section, we will show that, for such an $S$, $\e_{HK}(S)$ will turn out to be greater than $\e_{HK}(R)$ in many instances. 

We would like to outline this construction in a specific example.

Let $\ringR$ be the ring obtained by killing the $2 \times 3$ -minors of a generic matrix, say $R = k[[x,y,z,u,v,w]]/(xv-uy,yw-vz,xw-uz)$. This ring is Cohen-Macaulay of dimension $4$ with $x, u-y, z-v,w$ a system of parameters. In fact, $R$ is $F$-regular.

Let $A= k[[x, u-y, z-v, w]] \subset R$ be a Noether normalization. For computational purposes, let $a= u-y, b = z-v$. With this change of variables $A= k[[x,a,b,w]] \subset R = k[[x,a,b,w,y,v]]/(y^2-xv+ay,yw-vb -v^2, xw-ab-yv-av-yb)$. Note that $Q(A) \subset Q(B)$ is a simple field extension generated by $y$. Indeed, $v = \frac{1}{x} (y^2 + ay)$.

Look now at $A[[y]] \to R$. The kernel of this map is a principal ideal generated by some $f$. Hence we have constructed a hypersurface $\ringS$ in $R$. It is known that $\e(S) = \e(R)$. We would like to compare the Hilbert-Kunz multiplicities of $R$ and $S$. 

Since $R$ is finite over $S$, we have that $\e_{HK}(\fn, S) = \e_{HK}(\fn R, R)/ r$, where $r$ is the rank of $Q(R)$ over $Q(S)$ (by Theorem 2.7 in~\cite{WY1}). But $Q(S) = Q(R)$ and so $r=1$. We can also note that $\fn R \subset \fm $, which implies that $\e_{HK} (\fn R, R) \geq \e_{HK} (\fm, R) = \e_{HK}(R)$. Moreover, $R$ is $F$-regular and so $ \fn R = (\fn  R)^{*} \neq \fm$ which shows that $\e_{HK} (S) > \e_{HK}(R)$. ( As the referee points out, the reader can note that $\e_{HK}(R) = 13/8$ by applying the results of Section 5 in~\cite{WY3}.)

Examples like this are likely to abound. We have only used that $R$ is $F$-regular and that the finite extension $S \into R$ has rank $1$.

\section{Complete intersections}

In this section, we give an affirmative answer to the Conjecture~\ref{conjecture} i) in the case of complete intersections. We do this by reducing the study of complete intersections to that of hypersurfaces, a case that was solved in the previous section.

We would like to state first prime avoidance result that will be used later in the section (~\cite{Ei}, Exercise 3.19).

\begin{Lemma}[Prime Avoidance]
\label{prime}
Suppose that $R$ is a ring containing a field $k$, and let $I_1,...,I_m$ be ideals.
If $f_1,...,f_n \in R$ are such that $(f_1,...,f_n) \nsubseteq I_i$
 for each $i$,
 then there exists a nonzero homogeneous polynomial $H(Z_1,...Z_n) \in k[Z_1,...,Z_n]$ such that 
$$\sum_{i=1}^{n}a_if_i \notin \bigcup_i I_{i}$$ 
for all $(a_1,...,a_n) \in k^n$ with $H(a_1,...,a_n) \ne 0$.
\end{Lemma}

The Lemma will be used in the proof of the following

\begin{Proposition}
\label{sc-ci}
Let $k$ be an uncountable algebraically closed field of characteristic $p >0$. Let $A = k[[X_1,...,X_n]]$ and $\tilde R : = A/(f_1...f_l)$ a complete intersection ring and $f,g \in A$ such that they form a regular sequence on $\tilde R$. Let $0 \neq h \in \tilde R$.
Then there exist a dense subset $V \subset k $ such that $ah +f, g$ form a regular sequence on $\tilde R$ and
$$ \e_{HK} (\tilde R/ (f,g) ) \geq \e_{HK} (\tilde R/ (a h+f, g),$$
for all $a \in V$.
\end{Proposition}

\begin{proof}
Since $f,g$ form a regular sequence on $\tilde R$, we note that $(h,f) \not\subseteq P$ for every associated prime $P$ of $\tilde R/ (g)$. Hence, we can find a nonzero homogeneous polynomial $H(Z_1,Z_2)$ such that $$a h +f \notin P$$
for every $P$ associated prime of $\tilde R/ (g)$ and every $a$ in the open non-empty subset $U: = \{ a \in k: H(a,1) \neq 0 \}$. That is,  $ah +f$ and $ g$ form a regular sequence on $\tilde R$. Let us consider the natural ring homomorphism
$$k [t] \to \tilde R [t] / (th+f, g).$$
The fiber over each $a \in U$ is of dimension $n-l-2$. As in the proof of 
Theorem~\ref{sc-hyp} we can find a dense subset $V$ in $U$ such that

$$\e_{HK} (\tilde R / (f, g) \geq \e_{HK} (\tilde R/ (ah+f, g),$$
for all $a \in V$.
\end{proof}

\begin{Theorem} 
\label{ci}
Let $(R,\fm,k)$ be a non-regular complete intersection whose residue field is
an uncountable algebraically closed field of characteristic $p>0$.
Then there exists a non-regular hypersurface $k[[X_1,...,X_{d+1}]]/(F)$ such that
$$\e_{HK}(k[[X_1,...,X_{d+1}]]/(F)) \le \e_{HK}(R).$$ 
\end{Theorem}
\vspace{0.3cm}

\begin{proof}
Let $R$ be a non-regular complete intersection of dimension d.
Since we can complete $R$, $R$ is isomorphic to 
$$k[[X_1,...,X_{d+e}]]/(f_1,...,f_e),$$ 
where $(f_1,...,f_e)$ is a regular sequence.

\vspace{0.2cm}

($e=1$): In this case, since $R$ is already a hypersurface, so we are done.
\vspace{0.2cm}

($e>1$): We will give a proof based on induction on the length of a regular sequence.
The idea of the proof is to work on the regular sequence. In each step, we try to obtain
another regular sequence whose corresponding residue ring is of dimension d, non-regular, and has multiplicity
smaller than equal to that of the residue ring corresponding to regular sequence obtained in the previous step.

First of all,
we will apply the following procedures to the ring $R$.
\vspace{0.2cm}

(1): Suppose that some $f_i$ ($1 \le i \le e$) defines a regular hypersurface ring, then by Cohen's structure theorem, there is an isomorphism $$k[[Y_1,...,Y_{d+e-1}]] \cong k[[X_1,...,X_{d+e}]]/(f_i),$$ where $k[[Y_1,...,Y_{d+e-1}]]$ is the power series ring.
Then there is an isomorphism $$k[[Y_1,...,Y_{d+e-1}]]/(f'_1,...,f'_{i-1},f'_{i+1},...,f'_e) \cong k[[X_1,...,X_{d+e}]]/(f_1,...,f_e)),$$ 
where $f'_j$ is the inverse image of $f_j$. Note that $(f'_1,...,f'_{i-1},f'_{i+1},...,f'_e)$ is a regular ideal and its length is equal to $e-1$.  

Following this procedure, we can shrink the length of the regular sequence as small as possible, therefore we can assume that none of $f_i$'s defines a regular hypersurface.
\vspace{0.2cm}

(2): After (1) is done, by making some linear change of $X_1,...,X_{d+e}$, we can assume that each $f_i$ contains a term, $c_iX_{1}^{t_i}$ with $0\ne c_i\in k$, and that the order of $f_i$ is equal to $t_i$ for each $i$. 
The coefficients of $X_1^{t_i}$ are of the form $c_i+m_i$ with $m_i$ in the maximal ideal of $k[[Y_2,...,Y_{d+e-1}]]$.
Then by Weierstrass preparation theorem, each $f_i$ can be written uniquely in the form  $$f_i=u_i(X_1^{t_i}+a_{s-1}X_1^{t_i-1}+\cdots+a_0),$$ where $u_i$ is a unit, and $a_i$ is in the maximal ideal of $k[[Y_2,...,Y_{d+e-1}]]$.
\vspace{0.3cm}

Since we consider ideals, so we can ignore the unit $u_i$,
 hence again, we may put $$f_i=(X_1^{t_i}+a_{s-1}X_1^{t_i-1}+\cdots+a_0),  R:=k[[X_1,...,X_{d+e}]]/(f_1,...,f_{e}).$$
 To apply the induction step, let us prove the following proposition.

\begin{Proposition}
\label{udsc}
Let $\tilde R:=k[[X_1,...,X_n]]/(f_1,...,f_l)$ be a complete intersection and $f$, $g$ be elements of $A:=k[[X_1,...,X_n]]$ that form a regular sequence on $\tilde R$.
Assume that both $A/(f)$ and $A/(g)$ are non-regular, and $f$, $g$ are distinguished polynomials with respect $X_1$, that is, they can be written as $f=(X_1^{t}+a_{t-1}X_1^{t-1}+\cdots+a_0)$, $g=(X_1^{s}+b_{s-1}X_1^{s-1}+\cdots+b_0)$,
 where $a_i$, $b_i$ are in the maximal ideal of $k[[X_2,..,X_n]]$.

Then, there exists a regular sequence  $f', g'\in k[[X_1,...,X_n]]$ in $\tilde R$ such that 
$$ \e_{HK}(\tilde R/(f,g)) \ge \e_{HK}(\tilde R/(f',g')),$$
and such that following holds: 

$f'$ $(or ~g')$ contains a linear term in $X_1$: that is,  $f' = u' X_1 + v'$ with
$u'$ unit in $\tilde{R}$ and $v' \in k[[X_2,..., X_n]]$ 
 
Moreover, one can arrange that $\tilde R/(f', g')$ is  non-regular.
\end{Proposition}

\begin{Remark}
By Kunz, Proposition~\ref{kunz}, we note that $e_{HK}(\tilde R/(f)), e_{HK}(\tilde R/(g)) \le e_{HK}(\tilde R/(f, g))$, hence $\tilde R/(f, g)$ is also non regular. In the same manner, if one of $f'$ and $g'$ defines a non-regular hypersurface, then $\tilde R/(f', g')$ is also non-regular.
\end{Remark}

\begin{proof}[Proof of the Proposition] 

The plan is to start with the ideal $(f,g)$ in $\tilde{R}$ and perform transformations on $f$ or $g$ to decrease the degree of $X_1$ in either $f$ or $g$ until we come to one of the cases described below.

The first step is natural and easy to describe: Without loss of generality, we may assume $t \ge s$. 
Then $F':=f-X_1^{t-s}g$ has $deg_{X_1}(F') < t$, where $deg_{X_1}$ denotes the degree with respect to $X_1$.
So we have $(f,g)=(F',g)$ as ideals.
Since every $a_i$ and $b_i$ is in the maximal ideal, the top coefficient of $F'$ is also in the maximal ideal.
We see that $F', g$ is a regular sequence by the vanishing of Koszul homology.
Let us put $t':=deg_{X_1}(F')$, $s':=deg_{X_1}(g)$, and $G':=g$. So starting with $f, g$, we obtained $F', G'$.

This first step fits under the general procedure that is described in the next:

We have two elements  $F, G \in k[[X_1,...,X_n]]$ in $\tilde R$ such that 
$$ \e_{HK}(\tilde R/(f,g)) \ge \e_{HK}(\tilde R/(F,G)),$$ and, at least one of them, say $F$, has the leading term in $X_1$ of the form $u X_1^{s}$, with $u$ a unit in $\tilde{R}$.   

We would like to show that one can construct $F', G'$ such that 
$$ \e_{HK}(\tilde R/(F,G)) \ge \e_{HK}(\tilde R/(F',G')),$$ and 
 $deg_{X_1}(F)+deg_{X_1}(G) > deg_{X_1}(F')+deg_{X_1}(G')$, such that either $F'$ (or $G'$) has the leading term in $X_1$ of the form $u' X_1^{t'}$ (or $u' X_1^{s'}$) with $u'$ a unit.

The first step described above is a particular case of the general procedure if one takes $F:=f, G:=g$.

Let us explain now how to make $F', G'$ from the given $F, G$. Let $deg_{X_1}(F) =t $ and $deg_{X_1}(G) =s$ and as above $F = u X_1^t + \cdots$, with $u$ a unit in $\tilde{R}$ and $G = vX_1^s +\cdots$, with $v$ not necessarily a unit. 

We have two cases to consider for the ideal $(F,G)$ as follows.
\vspace{0.2cm}

($\alpha$): If $t \le s$, we can take $$G':=G-vX_1^{s-t}{u}^{-1}F,~F':=F,$$ and put $t':=deg_{X_1}(F')$, and $s':=deg_{X_1}(G')$. Then we see that $deg_{X_1}(G)>deg_{X_1}(G')$ and that $(F',G')=(F,G)$.
Again $F', G'$ is a regular sequence on $\tilde R$. 
\vspace{0.2cm}

($\beta$): If $t \ge s$, then we can not use $G$ to eliminate the leading term in $X_1$ in $F$ since $v$ might not be a unit. Hence we will use Proposition~\ref{sc-ci} to replace $G$ by another  power series $G_1$ such that $G_1$ has the leading term in $X_1$ of the form $v_1 X_1^s$ where $v_1$ is a unit in $\tilde{R}$.

Consider the sequence ${a}X_1^{s}+G$, $F$, where $a \in k$.
Note that the top coefficient of $G_1: = {a}X_1^{s_1}+G$ is a unit in $A$ unless
$a=0$.

We apply Proposition~\ref{sc-ci} for $A$, $\tilde R$ and the regular sequence $F, G$ on $\tilde R$: there is a dense subset $V \subseteq \Max(k[t]) \simeq k$ for which
$$\e_{HK}(\tilde R/(F, G)) \ge \e_{HK}(\tilde R/(aX_1^{s}+G, F))$$ holds for all $a \in V$, and $aX_1^{s}+G, F$ form a regular sequence.

Working with the new sequence $(F, G_1=a X_1^{t}+G)$ for some $a\ne 0$ and $a \in V$, we obtain a new regular sequence $F', G'$ such that

$$ F': = F - u X_1^{t-s}v_1^{-1}G_1, \ G' :=G_1 $$ where $v_1$ is the top coefficient of $G_1$. Also we remark that $(F', G')=(F, G_1)$ as ideals, and $deg_{X_1}(F)>deg_{X_1}(F').$

One can see in either case $F'$ (or $G'$) has the leading term in $X_1$ of the form $u' X_1^{t'}$ (or $u' X_1^{s'}$) with $u'$ a unit.

Moreover, the new pair $F',G'$ satisfies the property: $deg_{X_1}(F')+deg_{X_1}(G') < deg_{X_1}(F)+deg_{X_1}(G)$.  We also note that whenever we apply Proposition~\ref{sc-ci}, then the ideal $(F', G')$ is different than the ideal $(F, G)$.

Once we have $F', G'$, we continue by applying the procedure to $F', G'$ themselves.
We would like to show that by doing this repeatedly we will eventually reach one of the forms stated in the conclusion of the Proposition. 

Both $f,g$ belong to $\fm_A^2$. We notice that if $F, G$ belong to $\fm_A ^2$, then $F', G'$ will also belong to $\fm_A ^2$ unless $min(deg_{X_1}(F), deg_{X_1}(G))=1$. Once this situation occurs, we stop our procedure at once; if say $deg_{X_1}(F) =1$, then by changing the coefficient of $X_1$ with the help of Proposition~\ref{sc-ci} if necessary, we see that we end up in the case described. 

If we never encounter the situation where $min(deg_{X_1}(F), deg_{X_1}(G))=1$, then we eventually end up with $f'$ (or $g'$) $\in k[[X_2,...,X_n]].$  But then using Proposition~\ref{sc-ci} add $uX_1$ to $f'$ or $g'$ and we end up in the situation described in the conclusion
of our Proposition.

To end the proof, it is enough to say that at least one of $f'$ or $ g'$ is in $\fm_A^2$. Then this guarantees that $\tilde{R}/(f',g')$ is non-regular.
\end{proof}

\vspace{0.8cm}

Now let us go back to the proof of the theorem.
We apply the Proposition 4.4 for $A:=k[[X_1,...,X_{d+e}]]$, $l:=e-2$ to $f_1,...,f_{e}$
inductively.

Start with $f_1$ and $f_2$ and put $\tilde R:=k[[X_1,...,X_{d+e}]]/(f_3,...,f_{e})$.
Then we can find such $F_1, F_2$ as stated in the Proposition.
Once we come to the conclusion in the Proposition, then we can find the desired hypersurface by applying the induction step on the length of the regular sequence by eliminating $X_1$, so we are done.

\end{proof}

We would like to close this section by proving the part (1) of Conjecture of Watanabe and Yoshida stated in the introduction for complete intersections

\begin{Theorem}
\label{main}
Let $d \geq 2$, $ p \neq 2$ prime and $k$ a field of characteristic $p>0$. If $(R, \fm, k)$ is a complete intersection, not regular, then $\e_{HK}(R) \geq \e_{HK}(R_{d,p})$.
\end{Theorem}

\begin{proof}

We can enlarge the residue field such that we have an uncountable algebraic closed field $K$. 

By Theorems~\ref{hypersurface} and~\ref{ci} we see that over $K$, $\e_{HK} (R \otimes_k K) \geq \e_{HK} (R_{d,p} \otimes_k K)$ which implies the result over $k$.

\end{proof}

\begin{Remark}
{\rm Although we stated Propositions~\ref{udsc} and~\ref{sc-ci} for the case of complete intersection only, this assumption was in fact not needed in their corresponding proofs. We kept this as hypothesis for the convenience of the reader, since this section deals only with complete intersections.}

\end{Remark}

\section{Remarks on the general case}

In this section, we would like to show how using ideas related to the upper semi-continuity of the Hilbert-Kunz multiplicity can provide insight into the general case of the Conjecture stated in Section 1. A local ring $S$ such that $\dim (S) - {\rm depth} (S) =1 $ is called $almost ~Cohen$-$Macaulay$.

\begin{Proposition}
Let $\ringR$ be an catenary unmixed non-regular ring of positive characteristic $p >0$. Then there exists a non-regular unmixed ring of same dimension $\ringS$ which is Cohen-Macaulay or almost Cohen-Macaulay such that $$\e_{HK} (S) \leq \e_{HK} (R).$$
\end{Proposition}

\begin{proof}

Let $x_1, \cdots, x_n$ be a maximal regular sequence on $R$ and let $P$ be a minimal prime over $(x_1, \cdots, x_n)$. We have that $\e_{HK} (R_P) \leq \e_{HK}(R)$ by Theorem 3.8 in~\cite{K2} (this is where we need catenary). If $R_P$ is not regular we are done, since we can adjoin a finite number of indeterminates to $R_P$ to obtain a Cohen-Macaulay ring $S$ with $\e_{HK} (S) = \e_{HK} (R_P) \leq \e_{HK} (R)$ (the first equality comes from Proposition~\ref{kunz}).

If $R_P$ is regular, then consider $P \subset Q$ such that $\height(Q/P)=1$. Localize at $Q$ and get $\e_{HK}(R_Q) \leq \e_{HK} (R)$. Since $x_1, \cdots, x_n$ is a maximal regular sequence we see that $R_Q$ is almost Cohen-Macaulay. As before, by adjoining a number of indeterminates over $R_Q$ we obtain an example of same dimension as $R$.
\end{proof}

We would like to show that part (1) of the Conjecture can be reduced to the case of an isolated singularity:

Assume that $\ringR$ is excellent and unmixed. It is immediate that $\e_{HK}(R) \geq \e_{HK}(R_{red})$ and hence we can pass to $R_{red}$ and assume that $R$ is excellent and reduced.

By induction on the dimension of $R$ we can assume that for all non-regular unmixed rings $A$ of smaller dimension one can find a hypersurface $B$ of same dimension such that $\e_{HK} (B) \leq \e_{HK}(A)$.

Let $\Sing(R)$ be the singular locus of $\ringR$. It is a non-empty closed set defined by an ideal $J$. If $J$ is $\fm$-primary, then there is nothing to prove. Otherwise, let $P_i$, $i=1, \cdots, n$, be the collection of all minimal primes of $J$. Let $P$ be one such minimal prime $P_i$ with height less than the dimension of $R$.

Then $\e_{HK} (R_{P}) \leq \e_{HK} (R)$. By induction, we can find a hypersurface $S$ such that $\e_{HK} (S) \leq \e_{HK} (R_P)$. By adjoining a finite number of indeterminate to $R_P$ we obtain a hypersurface, relabeled $S$, of dimension equal to $\dim (R)$ and $\e_{HK}(S) \leq \e_{HK} (R)$.

Our result Theorem~\ref{hypersurface} shows that among hypersurfaces $\sum_{i=0}^{d} x_i ^2$ is the one with minimal Hilbert-Kunz multiplicity.

We would like to close now with an observation related to the questions addressed in this paper: Let $A$ be a finitely generated $K$-algebra which is non-regular and locally unmixed. Is there a minimal value for the Hilbert-Kunz multiplicity of $A_P$ where $P$ is a non-regular prime?

\begin{Proposition}

Let $A$ be an excellent, nonregular and locally unmixed. Then $\e_{HK} : \Spec (R) \to \mathbf{R}$ has minimum when restricted to the non-regular locus of $\Spec (R)$.

\end{Proposition}

\begin{proof}
$A$ is excellent and hence its singular locus is defined by an ideal $J$. For any prime containing $J$ we can find a minimal prime $P$ of $J$, $P \subset Q$ such that $\e_{HK} (A_P) \leq \e_{HK} (A_Q)$. 

Since there are only finitely many minimal primes over $J$ we are done.
\end{proof}

\end{document}